\documentclass{amsart}

\usepackage{amssymb}
\usepackage{amsthm}
\usepackage{amsmath}
\usepackage{enumitem}
\usepackage[mathscr]{euscript}
\usepackage{mathtools}
\usepackage{stmaryrd}
\usepackage{cancel}
\usepackage{soul}
\usepackage{tikz}
\usepackage{tikz-cd}
\usepackage[normalem]{ulem}
\usepackage[unicode]{hyperref}
\hypersetup{colorlinks, linkcolor=blue}
\usepackage{cleveref}
\usepackage{tikz}
\usetikzlibrary{fadings}
\usepackage{float}
\usetikzlibrary{patterns}
\usetikzlibrary{shadows.blur}
\usetikzlibrary{shapes}
\usepackage{xcolor}
\usepackage{todonotes}
\usepackage{comment}

\theoremstyle{plain}

\newtheorem{proposition}{Proposition}[section]
\newtheorem{theorem}[proposition]{Theorem}
\newtheorem{lemma}[proposition]{Lemma}

\theoremstyle{definition}

\newtheorem{remark}[proposition]{Remark}

\DeclareMathOperator{\SL}{\mathsf{SL}}

\DeclareMathOperator{\Fc}{\mathcal{F}}

\DeclareMathOperator{\Mmc}{\mathcal{M}}

\DeclareMathOperator{\Pc}{\mathcal{P}}

\DeclareMathOperator{\Cb}{\mathbb{C}}

\DeclareMathOperator{\Nb}{\mathbb{N}}

\DeclareMathOperator{\Rb}{\mathbb{R}}

\DeclareMathOperator{\Zb}{\mathbb{Z}}

\newcommand{\abs}[1]{\left|#1\right|}

\begin{document}

\title[A non-Morse relatively Anosov group]{A non-Morse relatively Anosov group}

\author{Tianqi Wang}
\address{Department of Mathematics, Yale University, New Haven, CT, USA}
\email{tq.wang@yale.edu}
\thanks{}

\date{\today}

\begin{abstract}
    We construct a relatively Anosov group that is not relatively Morse.
\end{abstract}

\maketitle


\section{Introduction}

The notion of Anosov group, introduced by Labourie \cite{labourie2004anosov} and later developed in \cite{guichard2012anosov,guichard2017anosov,kapovich2014morse,kapovich2017anosov,bochi2019anosov}
serves as a natural generalization of convex-compact group in higher rank semisimple Lie groups. Given a finitely generated, $\mathsf{P}_1$-divergent, discrete subgroup $\Gamma< \SL(d,\Cb)$, we denote by $\Lambda(\Gamma)$ its limit set in $\Fc_{1,d-1}$ of $(1,d-1)$-partial flags. 
Kapovich--Leeb--Porti \cite{kapovich2017anosov} characterized the $\mathsf P_1$-Anosov property of $\Gamma$ by either of the following equivalent conditions, \begin{itemize}
    \item $\Lambda(\Gamma)$ is transverse and the action of $\Gamma$ on $\Lambda(\Gamma)$ is a uniform convergence action;
    \item an orbit map from $\Gamma$, equipped with a word metric, to the symmetric space of $\SL(d,\Cb)$ is a $\mathsf P_1$-Morse quasi-isometric embedding.
\end{itemize}

These two approaches also lead to generalizations of geometric finiteness to higher-rank semisimple Lie groups. We retain the notation above. From the point of view of boundary actions, $\Gamma$ is relatively $\mathsf P_1$-asymptotically embedded if $\Lambda(\Gamma)$ is transverse and the action of $\Gamma$ on $\Lambda(\Gamma)$ is a geometrically finite convergence action in the sense of Bowditch \cite{bowditch2012relatively} and Yaman \cite{yaman2004}. This notion was introduced by Kapovich--Leeb \cite{kapovich2023relativizing} and is also referred to as relative $\mathsf P_1$-Anosovness as studied in Zhu--Zimmer \cite{zhu2022relatively}. From the point of view of orbit maps, Kapovich--Leeb \cite{kapovich2023relativizing} defined $\Gamma$ to be relatively $\mathsf P_1$-Morse if there exists $(\Gamma,\Pc)$, a relatively hyperbolic structure on $\Gamma$, a Gromov model $Y$ of $(\Gamma,\Pc)$, and a $\Gamma$-equivariant $\mathsf P_1$-Morse quasi-isometric embedding of $Y$ into the symmetric space of $\SL(d,\Cb)$. We use the Gromov model convention of \cite{kapovich2023relativizing} and when refer to a \emph{roughly geodesic} Gromov model, we allow uniform rough geodesics in place of geodesics.

Kapovich--Leeb \cite[Theorem~8.3]{kapovich2023relativizing} proved that every relatively $\mathsf P_1$-Morse subgroup is relatively $\mathsf P_1$-Anosov. Zhu--Zimmer introduced uniform relative Anosov subgroups with respect to a Gromov model $X$ \cite[Definition~1.11]{zhu2022relatively}. They proved that a uniform relative $\mathsf P_1$-Anosov group with respect to $X$ determines a $\Gamma$-equivariant $\mathsf P_1$-Morse quasi-isometric embedding of $X$ into the symmetric space \cite[Proposition~11.3]{zhu2022relatively}. Conversely, if a discrete subgroup admits such an embedding that is continuous, it is uniform relative $\mathsf P_1$-Anosov \cite[Proposition~11.4]{zhu2022relatively}.

It is therefore natural to ask whether every relatively $\mathsf P_1$-Anosov subgroup is relatively $\mathsf P_1$-Morse. In this paper, we give a negative answer by constructing a counterexample using the following facts.

\subsection*{A ping-pong argument}

Let $P<\SL(d,\Cb)$ be a finitely generated, discrete, $\mathsf P_1$-divergent, weakly unipotent subgroup whose limit set consists of a single flag $F$. Choose a $\mathsf P_1$-biproximal element $g\in\SL(d,\Cb)$ whose attracting and repelling flags are transverse to $F$. By Zhu--Zimmer \cite[Proposition~12.2]{zhu2024relatively}, there exists a subgroup $P' < P$ of finite index such that for $N\in\Nb$ large enough, the subgroup generated by $g^N$ and $P'$ is a free product \[\Gamma=\langle g^N\rangle * P'\] that is relatively $\mathsf P_1$-Anosov relative to $\{P'\}$. In particular, $(\Gamma,\{P'\})$ is relatively hyperbolic.

\subsection*{Hyperbolic scalar Cartan metrics}

In joint work with Reyes \cite[Theorem~D]{reyes2026regularity}, we proved that for relatively Morse groups, the scalar Cartan metrics are hyperbolic and give rise to roughly geodesic Gromov models. We interpret the theorem in our setting as below.

\begin{lemma}\label{lemma: hyperbolicity}
    If $\Gamma<\SL(d,\Cb)$ is $\mathsf P_1$-Morse relative to $P$, then $(\Gamma,P)$ admits a roughly geodesic Gromov model $(Y,d_Y)$, a basepoint $o\in Y$, and a constant $C_1\geqslant0$ such that for any $\gamma \in \Gamma$,
\[\abs{d_Y(o,\gamma o)-\log\frac{\sigma_1(\gamma)}{\sigma_d(\gamma)}} \leqslant C_1~.\]
\end{lemma}

Here $\sigma_1(\gamma)\geqslant\cdots\geqslant\sigma_d(\gamma)>0$ denote the singular values of $\gamma\in\SL(d,\Cb)$.

\subsection*{Almost ultrametric peripheral subgroups}
A metric space $(H,d_H)$ is said to be \emph{almost ultrametric} if there exists a constant $C\geqslant 0$ such that for any $x,y,z \in H$, \[d_H(x,z)\leqslant\max\{d_H(x,y),d_H(y,z)\}+C~.\] we use a lemma from unpublished joint work with Reyes, which states any horosphere in a roughly geodesic Gromov model is almost ultrametric. More concretely, 

\begin{lemma}\label{lemma: ultrametric}
    Suppose $(Y,d_Y)$ is a roughly geodesic Gromov model of $(\Gamma,P)$ and $o\in Y$ is a fixed basepoint. Then there exists a constant $C_2\geqslant 0$ such that for any $\gamma,\eta\in P$, \[d_Y(o, \gamma \eta o ) \leqslant \max\{d_Y(o, \gamma o ),d_Y( \gamma o, \gamma \eta o ) \}+C_2\]
\end{lemma}

The lemma follows from a standard hyperbolic geometry argument (see \cite[Proposition~2.6]{wang2023notions}) and we include a proof in Section~\ref{section: ultrametric}.\\

\subsection*{Construction of the example}

We construct a finitely generated, $\mathsf{P}_1$-divergent, weakly unipotent subgroup $P$ of $\SL(d,\Cb)$ with limit set a singleton such that \[\gamma \mapsto \left(\log\frac{\sigma_1(\gamma^2)}{\sigma_d(\gamma^2)}- \log\frac{\sigma_1(\gamma)}{\sigma_d(\gamma)}\right)\] is unbounded on $P$ in Theorem~\ref{theorem: the example}.

Then an example of relatively $\mathsf P_1$-Anosov groups that is not relatively $\mathsf P_1$-Morse can be construted in the following way.
Note that any finite index subgroup of $P$ also has these properties. Therefore by the ping-pong argument discussed above, we may assume $P$ is the peripheral subgroup of a relatively $\mathsf P_1$-Anosov group $\Gamma$. We claim that $\Gamma$ is not relatively $\mathsf P_1$-Morse. Suppose otherwise, then by Lemma~\ref{lemma: hyperbolicity} and Lemma~\ref{lemma: ultrametric}, for any $\gamma\in P$, \[\log\frac{\sigma_1(\gamma^2)}{\sigma_d(\gamma^2)}- \log\frac{\sigma_1(\gamma)}{\sigma_d(\gamma)}\leqslant 2C_1 +C_2~,\] where $C_1\geqslant 0$ and $C_2\geqslant 0$ are the constants given by the two lemmas respectively, which gives a contradiction.

\begin{theorem}\label{theorem: the example}
    There exists a discrete subgroup $P <\SL(5,\Cb)$ such that 
    \begin{enumerate}
        \item $P \cong \Zb^2$;
        \item $P$ is unipotent in $\SL(5,\Cb)$;
        \item $P$ is $\mathsf{P}_1$-divergent and $\Lambda(P)\subset \Fc_{1,4}$ is a singleton;
        \item There exists a sequence $(\gamma_n)_{n\in \Nb}\subset P$ such that \[\lim_{n\to \infty} \left(\log\frac{\sigma_1(\gamma_n^2)}{\sigma_5(\gamma_n^2)}- \log\frac{\sigma_1(\gamma_n)}{\sigma_5(\gamma_n)}\right) = +\infty~.\] 
    \end{enumerate}
\end{theorem}

\medskip

\subsection*{Acknowledgments}
The author would like to thank Balthazar Fl\'echelles for introducing him to the construction of Cooper's example, Eduardo Reyes for many interesting discussions on the almost ultrametric condition, and Andrew Zimmer for the discussion on several possible equivalent characterizations of the Morse property.

This paper is partially based upon work supported by the National Science Foundation under Grant No. DMS-2424139, while the author were in residence at the Simons Laufer Mathematical Sciences Institute in Berkeley, California, during the Spring 2026 semester.

\subsection*{Tool Use}
The mathematical ideas and arguments in this paper are developed, written, and checked by the author. A Large Language Model was used to assist with proofreading and correcting grammar and typos.

\section{Constructing the example}

Our construction is inspired by Cooper's example \cite{cooper2017} of a faithful unipotent representation of the Heisenberg group into $\SL(10,\Rb)$ that preserves a strictly convex domain in $\mathbb{RP}^9$, which can be viewed as a representation of translations on polynomials.

For $s,t\in\Cb$, let $T_{s,t} : \Cb[x,y] \to \Cb[x,y]$ be the translation defined by $(T_{s,t}f)(x,y) = f(x+s,y+t)$ for any $f\in \Cb[x,y]$. We fix the polynomial $f(x,y) = x^2 +y^3$. Then 
\[(T_{s,t}f)(x,y) = f(x+s,y+t) = f(x,y)+2sx+3ty^2+3t^2y+(s^2+t^3)~.\]

The subspace of $\Cb[x,y]$ spanned by $1,x,y,y^2,f$ is $T_{s,t}$-invariant. We identify this subspace with $\Cb^5$ with basis $e_1 = 1$, $e_2 = x$, $e_3 = y$, $e_4 = y^2$, $e_5 =f$. Then the restriction of $T_{s,t}$ on this subspace has matrix \[M_{s,t} = \begin{bmatrix} 1 & s & t & t^2 & s^2+t^3 \\ & 1 & 0 & 0 & 2s \\ & & 1 & 2t & 3t^2 \\ & & & 1 & 3t \\ &&&& 1\end{bmatrix}~.\]
Note that the $(1,2)$ and $(1,3)$ entries recover $s$ and $t$ respectively. Therefore $T_{s,t}\mapsto M_{s,t}$ is a faithful unipotent representation, \[\Mmc= \{M_{s,t}: s,t\in \Cb\}\cong \Cb^2 \cong \Rb^4\] is abelian and unipotent in $\SL(5,\Cb)$, and in particular for any $s,t,u,v\in \Cb$, \[M_{s,t}M_{u,v}=M_{s+u,t+v}~, \quad M_{s,t}^{-1}=M_{-s,-t}~.\]

We show the subgroup \[P = \langle M_{1,0}, M_{i,1}\rangle = \{M_{m+in,n} : m,n \in \Zb\}< \Mmc~,\] where $i = \sqrt{-1}$, satisfies the conditions in Theorem~\ref{theorem: the example}. Conditions (1) and (2) are clear. We show (3) and (4) by the next two lemmas.

\begin{lemma}\label{lemma: divergent}
    $P$ is $\mathsf{P}_1$-divergent and the limit set of $P$ in $\Fc_{1,4}$ is a singleton.
\end{lemma}

\begin{proof}
    We claim that \[\frac{1}{f(m+in,n)}M_{m+in,n}\to E_{15} \quad \text{as }\abs{m}+\abs{n}\to\infty~,\] where $E_{15}$ is the matrix with $(1,5)$-entry $1$ and all other entries $0$. Set $L_0=\langle e_1\rangle$ and $H_0=\langle e_1,e_2,e_3,e_4\rangle$. For any nonzero vector $v=\sum_{j=1}^5v_je_j$ such that $\langle v\rangle$ is transverse to $H_0$, we have $v_5 \ne 0$. Then the claim implies $M_{m+in,n}[v]\to[e_1]$ in $\mathbb{CP}^4$ as $\abs{m}+\abs{n}\to\infty$. Moreover, this convergence is uniform on any compact subset of $\mathbb{CP}^4\setminus\mathbb P(H_0)$. By a similar argument in the dual space, we have for any hyperplane $H$ that is transverse to $L_0$, $M_{m+in,n}H\to H_0$ and the convergence is uniform on any compact subset of hyperplanes transverse to $L_0$.
    Therefore, by a standard dynamic argument in flag manifolds (see for instance \cite[Lemma~2.2]{canary2022cusped} and \cite[Lemma~5.2]{wang2023notions}), $P$ is $\mathsf P_1$-divergent and $\Lambda(P) = \{(L_0,H_0)\}$.

    We prove the claim by showing the $(1,5)$-entry of $M_{m+in,n}$ dominates the other entries when $\abs{m}+\abs{n}\to \infty$. Suppose otherwise there exists a constant $C>0$ and an escaping sequence of integer pairs $(m,n)$ such that \[\frac{|(m+in)^2+n^3|}{1+|m+in|+n^2}\leqslant C~.\] If $\abs{n}$ is bounded up to subsequence, then $\abs{m}\to\infty$ up to further subsequence, and the quotient goes to $+\infty$, which gives a contradiction. Hence we may assume up to subsequence $\abs{n}\to\infty$. 
    By direct computation, the inequality gives that \[\abs{2mn}\leqslant \sqrt{ \abs{m^2-n^2+ n^3}^2+\abs{2mn}^2 } = \abs{(m+in)^2+n^3}\leqslant C(1+ \abs{m}+\abs{n}+n^2)~.\]Then there is a constant $C'>0$, such that for $n$ large enough, \[\abs{m}\leqslant\frac{C(1+\abs{n}+n^2)}{2\abs{n}-C}\leqslant C'\abs{n}~.\] Therefore for $n$ large enough, \[\abs{(m+in)^2+n^3}\geqslant \frac{1}{2}\abs{n^3} \quad \text{and} \quad \abs{m+in} \leqslant (1+C')\abs{n}~,\] and hence \[\frac{\abs{(m+in)^2+n^3}}{1+ \abs{m+in}+\abs{n^2}} \geqslant \frac{\abs{n^3}}{2(1+(1+C')\abs{n}+\abs{n^2})} \to +\infty \quad \text{as } \abs{n}\to\infty~,\] which is a contradiction.
\end{proof}

\begin{remark}
    The generators of $P$ are chosen so that the $(1,5)$-entries of some sequence of elements have ``proper cancellation'' to violate the almost ultrametric inequality, while keeping $\mathsf P_1$-divergence as seen in the proof of Lemma~\ref{lemma: divergent}.
    For comparison, one can show $\langle M_{1,0},M_{0,i}\rangle$ is almost ultrametric and $\langle M_{1,0},M_{0,1}\rangle$ is not $\mathsf P_1$-divergent.
\end{remark}

\begin{lemma}\label{lemma: non-ultrametric}
    Let $\gamma_n = M_{n^3+in^2,n^2}$. Then \[\lim_{n\to \infty} (\log\frac{\sigma_1(\gamma_n^2)}{\sigma_5(\gamma_n^2)}- \log\frac{\sigma_1(\gamma_n)}{\sigma_5(\gamma_n)}) = +\infty~.\] 
\end{lemma}

For any $g\in \SL(d,\Cb)$, we denote by $\mathrm{L}(g)$ the maximum of the norms of the entries of $g$. Since any two norms on the space of matrices are equivalent, there is a constant $C \geqslant 1$ such that
\begin{equation}\label{equation: singular value}
    C^{-1} \mathrm{L}(g) \leqslant \sigma_1(g) \leqslant C \mathrm{L}(g)\quad \text{and} \quad
    C^{-1} \mathrm{L}(g^{-1}) \leqslant \frac{1}{\sigma_d(g)} \leqslant C \mathrm{L}(g^{-1})~.
\end{equation}

\begin{proof}[Proof of Lemma~\ref{lemma: non-ultrametric}]
    From the proof of Lemma~\ref{lemma: divergent} we know that \[\mathrm{L}(\gamma_n) = \abs{(n^3+in^2)^2 + (n^2)^3} = \abs{2n^6+2in^5-n^4}= 2\abs{n^6}+O(n^5)\] for large enough $n$. Here $O(\cdot)$ refers to the asymptotic growth rate for $n\to\infty$ as in calculus. Note that $M_{s,t}^k = M_{ks,kt}$ for any $k\in \Zb$ and $s,t\in \Cb$. We have that \[\mathrm{L}(\gamma_n^{-1}) = \abs{(-n^3-in^2)^2 - (n^2)^3} = \abs{2in^5-n^4}= 2\abs{n^5 }+ O(n^4)~,\]
    \[\mathrm{L}(\gamma_n^2) = \abs{(2n^3+2in^2)^2+(2n^2)^3} = 12\abs{n^6}+O(n^5)~,\]
    \[\text{and} \quad \mathrm{L}((\gamma_n^2)^{-1}) = \abs{(-2n^3-2in^2)^2-(2n^2)^3} = 4\abs{n^6}+O(n^5)~.\]
    By \eqref{equation: singular value}, \[\frac{\sigma_1(\gamma_n)}{\sigma_5(\gamma_n)} \leqslant C^{2}L(\gamma_n)L(\gamma_n^{-1}) = C^{2}4\abs{n^{11}} + O(n^{10})\]
    \[\text{and} \quad \frac{\sigma_1(\gamma_n^2)}{\sigma_5(\gamma_n^2)} \geqslant C^{-2}L(\gamma_n^2)L((\gamma_n^2)^{-1}) = C^{-2}48\abs{n^{12}} + O(n^{11})~,\]
    Therefore \[\left(\log\frac{\sigma_1(\gamma_n^2)}{\sigma_5(\gamma_n^2)}- \log\frac{\sigma_1(\gamma_n)}{\sigma_5(\gamma_n)}\right) = \log n + O(1) \to +\infty \quad \text{as } n\to \infty~.\]
\end{proof}

\section{Parabolic orbit metrics are almost ultrametric}
\label{section: ultrametric}

Suppose $(\Gamma,P)$ is a relatively hyperbolic group and $(Y,d_Y)$ is a roughly geodesic Gromov model of $(\Gamma,P)$. Let $K\geqslant 0$ be a constant such that for any distinct $x,y\in Y\cup \partial_\infty Y$, there exists a $K$-rough geodesic (that is, a $(1,K)$-quasi-geodesic) joining $x$ and $y$. Let $\delta\geqslant 0$ be a constant such that any $K$-rough geodesic triangle or ideal triangle is $\delta$-thin. Let $o\in Y$ be a fixed basepoint, and denote by $p\in\partial_\infty Y$ the parabolic point fixed by $P$.

The existence of horofunctions on roughly geodesic hyperbolic spaces follows from standard Busemann estimates (see for instance \cite[Chapter~3]{buyalo2007elements}).

\begin{lemma}\label{lemma: horofunction}
    There exists a function $h:Y\to\Rb$ with $h(o)=0$ and a constant $C\geqslant0$ such that
    \begin{enumerate}
        \item $\abs{h(x)}\leqslant C$ for any $x\in Po$;
        \item $\abs{h(y)-h(z)}\leqslant d_Y(y,z)+C$ for any $y,z\in Y$;
        \item $\abs{h(z)-h(x)-d_Y(x,z)}\leqslant C$ for any $x\in Po$ and $z$ on a $K$-rough geodesic ray from $x$ to $p$.
    \end{enumerate}
\end{lemma}

We now prove Lemma~\ref{lemma: ultrametric}.

\begin{proof}[Proof of Lemma~\ref{lemma: ultrametric}]
    Let $\ell_{o,p}$ be a $K$-rough geodesic ray from $o$ to $p$ in $Y$. Then $\ell_{\gamma \eta o,p}=\gamma\eta\ell_{o,p}$ is a $K$-rough geodesic ray from $\gamma\eta o$ to $p$ in $Y$. Let $h$ be the horofunction centered at $p$ and $C\geqslant 0$ be the constant given by Lemma~\ref{lemma: horofunction}.
    Let $\ell_{o,\gamma\eta o}$ be a $K$-rough geodesic segment joining $o$ and $\gamma\eta o$. Let $m$ be a coarse midpoint of $\ell_{o,\gamma\eta o}$, that is, $m\in \ell_{o,\gamma\eta o}$ is such that \[\abs{d_Y(o,m) - \frac{1}{2} d_Y(o,\gamma \eta o)}\leqslant 3K \quad \text{and} \quad \abs{d_Y(m,\gamma \eta o) - \frac{1}{2} d_Y(o,\gamma \eta o)}\leqslant 3K~.\]

    Since the ideal triangle with edges $\ell_{o,p}$, $\ell_{o,\gamma \eta o}$ and $\ell_{\gamma \eta o,p}$ is $\delta$-thin, there exists $v\in \ell_{o,p} \cup \ell_{\gamma \eta o,p}$ such that $d_Y(m,v) \leqslant \delta$. Without loss of generality, we assume $v\in \ell_{o,p}$. Then by Lemma~\ref{lemma: horofunction} and triangle inequality,
    \begin{equation}\label{equation: midpoint low}
    \begin{array}{ll}
        h(m) & \geqslant h(v) - \delta - C \geqslant h(o) +d_Y(o,v) -\delta - 2C \geqslant d_Y(o,m) - 2\delta -2C \\ & \geqslant  \frac{1}{2} d_Y(o,\gamma \eta o) - 2\delta -2C - 3K~.
    \end{array}
    \end{equation}
    Let $\ell_{o,\gamma o}$ (respectively $\ell_{\gamma o,\gamma\eta o}$) be a $K$-rough geodesic joining $o$ and $\gamma o$ (respectively $\gamma o$ and $\gamma\eta o$).
    We then consider the triangle with edges $\ell_{o,\gamma o}$, $\ell_{\gamma o,\gamma\eta o}$ and $\ell_{o,\gamma \eta o}$. There exists a point $w\in \ell_{o,\gamma o}\cup \ell_{\gamma o,\gamma\eta o}$, such that $d_Y(m,w) \leqslant \delta$. Without loss of generality, we assume $w \in \ell_{o,\gamma o}$
    Then by Lemma~\ref{lemma: horofunction} and triangle inequality,
    \begin{equation}\label{equation: midpoint up}
    \begin{array}{ll}
        h(m)  & \leqslant  h(w) + \delta + C \\
        & \leqslant \min\{h(o)+d_Y(o,w)+C\ , h(\gamma o)+ d_Y(w, \gamma o)+C\} +\delta \\  
        & \leqslant \min\{d_Y(o,w)\ , d_Y(w, \gamma o)\} +\delta + 2C \\
        & \leqslant \frac{1}{2}d_Y(o,\gamma o) + 3K + \delta + 2C \\
        & \leqslant \frac{1}{2}\max\{d_Y(o,\gamma o)\ ,d_Y(\gamma o,\gamma \eta o)\} + \delta + 2C + 3K~.
    \end{array}
    \end{equation}
    Therefore \eqref{equation: midpoint low} and \eqref{equation: midpoint up} give
    \[ d_Y(o,\gamma \eta o) \leqslant \max\{d_Y(o,\gamma o)\ ,d_Y(\gamma o,\gamma \eta o)\} + 6\delta + 8C + 12K~,\] which completes the proof.
\end{proof}

\bibliographystyle{alpha}
\bibliography{reference.bib}

\end{document}